\documentclass[12pt]{article}
\usepackage{amsmath}
\usepackage{graphicx,psfrag,epsf}
\usepackage{enumerate}
\usepackage{natbib}
\usepackage{url} 

\usepackage{amsfonts}
\usepackage{amsthm}
\usepackage{cleveref}
\theoremstyle{plain}
\newtheorem{theorem}{Theorem}[section]
\newtheorem*{theorem*}{Theorem}

\newcommand{\blind}{0}

\begin{document}

\def\spacingset#1{\renewcommand{\baselinestretch}%
{#1}\small\normalsize} \spacingset{1}


\if0\blind
{
  \title{\bf A new tail bound for the sum of bounded independent random variables}
  \author{Jackson Loper\\
    Department of Statistics, University of Michigan\\
    and \\
    Jeffrey Regier \\
    Department of Statistics, University of Michigan}
  \maketitle
} \fi

\if1\blind
{
  \bigskip
  \bigskip
  \bigskip
  \begin{center}
    {\LARGE\bf A new tail bound for the sum of bounded independent random variables}
\end{center}
  \medskip
} \fi

\bigskip
\begin{abstract}
We construct a new tail bound for the sum of independent random variables for situations in which the expected value of the sum is known and each random variable lies within a specified interval, which may be different for each variable.  This new bound can be computed by solving a two-dimensional convex optimization problem.  Simulations demonstrate that the new bound is often substantially tighter than Hoeffding's inequality for cases in which both bounds are applicable.
\end{abstract}

\noindent%
\vfill

\newpage
\spacingset{1.45} 

\section{Introduction}

The celebrated Hoeffding's inequality \citep[Theorem 2]{hoeffding_probability_1963} has been a cornerstone for probability and statistics over the last 60 years.  This inequality, which we refer to as the ``general Hoeffding inequality,'' concerns the sum of $n$ independent random variables, $(X_1,\ldots, X_n)$.  Each variable $X_i$ lies in some known interval $[a_i,b_i]$, for $i\in \{1,\ldots n\}$.  In this article, without loss of generality, we assume that $a_i=0$ for each $i$.  The general Hoeffding inequality states that the sum of these random variables, $S=\sum_{i=1}^n X_i$, is unlikely to greatly exceed its expected value, $\mu=\mathbb{E}[S]$.  In particular, for any $t\geq 0$,
\begin{equation} \label{eq:hoeffgeneral}
\mathbb{P}(S \geq \mu+t) \leq \exp\left(-2t^2 / \sum_{i=1}^n b_i^2 \right).
\end{equation}
Hoeffding's seminal work also included a less well-known bound, which is applicable only if every bounding interval is of the same length and the expected value $\mu$ is known \citep[Theorem 1]{hoeffding_probability_1963}.  If $b_i=1$ for every $i$, then
\begin{equation} \label{eq:hoeffspecial}
\mathbb{P}(S \geq \mu+t) \leq \left(\frac{\mu}{\mu+t}\right)^{n\mu+nt} \left(\frac{1-\mu}{1-\mu-t}\right)^{n-n\mu-nt}.
\end{equation}
We refer to \Cref{eq:hoeffspecial} as the ``specialized Hoeffding inequality.'' \Cref{eq:hoeffspecial} is tighter than \Cref{eq:hoeffgeneral} in situations where both are applicable.  

In this article, we propose a computable bound for the sum of independent random variables for situations in which the expected value of the sum is known and each random variable is known to lie in a fixed interval of arbitrary size, which may differ among variables.  Such bounds are useful for hypothesis testing regarding the sum's expected value based on an observed realization.  The new bound can be calculated by solving a two-dimensional convex optimization problem.  The result reduces to \Cref{eq:hoeffspecial} in the special case that $b_i=1$ for each $i \in \{1,\ldots,n\}$.  We use simulations to demonstrate that the new bound can be much tighter than the general Hoeffding inequality in situations where both are applicable and the specialized Hoeffding inequality is not.

\section{Notation and preliminaries}

\label{sec:notation}

Consider a set of probability distributions on independent bounded random variables.  For any scalar $c\geq 0$, let $\mathcal{M}_c$ denote the space of probability distributions with support on $[0,c]$.  For any vector of interval lengths $\boldsymbol{b}\in ({\mathbb{R}^+})^m$ and fixed value $\mu$, let
\begin{equation} \label{eq:mdef}
\mathcal{M}_{\boldsymbol{b},\mu} = \left\{\boldsymbol{p} \in \prod_{i=1}^m  \mathcal{M}_{b_i}:\ \sum_{i=1}^m \mathbb{E}_{X_i \sim p_i}[X_i] \leq \mu\right\}.
\end{equation}
In other words, a random vector $\boldsymbol{X}=(X_1,\ldots X_n)$  drawn according to a member of $\mathcal{M}_{\boldsymbol{b},\mu}$ satisfies three conditions: (i) each variable is independent, (ii) $X_i \in [0,b_i]$ for each $i$, and (iii) the expected sum over the vector is at most $\mu$.

\citet[Theorem 2]{hoeffding_probability_1963} provides a tail bound for the variable $S=\sum_{i=1}^m X_i$ when $\boldsymbol{X}$ is drawn according to $\mathcal{M}_{\boldsymbol{b},\mu}$, for fixed $\boldsymbol{b}$ and $\mu$.  To construct this bound, \citet{hoeffding_probability_1963} first considers a tail bound for a particular $\boldsymbol{p} \in \mathcal{M}_{\boldsymbol{b},\mu}$.
Let
\[
\varphi(\boldsymbol{p},s) = \inf_{t\geq 0} \left(\sum_{i=1}^m \log \mathbb{E}_{X_i \sim p_i}[\exp(t X_i)] - ts\right).
\]
If $\boldsymbol{X}$ is drawn according to any member of $\mathcal{M}_{\boldsymbol{b},\mu}$, the Chernoff--Cram\'er bound gives that $\log \mathbb{P} \left (S \geq s\right)$ is at most $\sup_{\boldsymbol{p} \in \mathcal{M}_{\boldsymbol{b},\mu}} \varphi(\boldsymbol{p},s)$.
To produce a bound that is universal for all members of $\mathcal{M}_{\boldsymbol{b},\mu}$, \citet{hoeffding_probability_1963} therefore considers the quantity
\begin{equation} \label{eq:varphistar}
\varphi^*_{\boldsymbol{b},\mu}(s) = \sup_{\boldsymbol{p} \in \mathcal{M}_{\boldsymbol{b},\mu}} \varphi(\boldsymbol{p},s) = \sup_{\boldsymbol{p} \in \mathcal{M}_{\boldsymbol{b},\mu}} \inf_{t\geq 0} \left(\sum_{i=1}^m \log \mathbb{E}_{X_i \sim p_i}[\exp(t X_i)] - ts\right).
\end{equation}
This leads to a tail bound for $S$ in terms of $\varphi^*$: 
\begin{equation} \label{eq:chernoff}
\mathbb{P} \left (S \geq s\right) \leq \exp(\varphi^*_{\boldsymbol{b},\mu}(s)).
\end{equation}
\Cref{eq:chernoff} is the worst-case Chernoff--Cram\'er bound for a sum of independent random variables $(X_1,\ldots X_n)$ with each $X_i \in [0,b_i]$ and $\sum_{i=1}^n \mathbb{E}[X_i]=\mu$.  

To compute the right-hand side of \Cref{eq:chernoff}, \citet{hoeffding_probability_1963} considers two different readily computable expressions for $\varphi^*$.  In one theorem, a formula for $\varphi^*_{\boldsymbol{b},\mu}$ is provided for the special case in which $b_1=\cdots=b_m=1$, leading to the specialized Hoeffding inequality, \Cref{eq:hoeffspecial}.  In another theorem, an \emph{upper bound} is produced for $\varphi^*_{\boldsymbol{b},\mu}$ in the general case of arbitrary $\boldsymbol{b}$, namely,
$
\varphi^*_{\boldsymbol{b},\mu}(s)\leq -2(s-\mu)^2/\sum_{i=1}^m b_i^2,
$
leading to the general Hoeffding inequality.  

Our contribution in the next section is to show how $\varphi^*$ can be \textit{exactly} computed for arbitrary interval lengths $\boldsymbol{b}$.
By substituting this exact value of $\varphi^*$ into \Cref{eq:chernoff}, an improved tail bound for the sum of bounded independent random variables follows.

\section{Optimizing the worst-case Chernoff--Cram\'er bound}

We now show that $\varphi^*$ (cf. Equation~\ref{eq:varphistar}) can be computed by solving a two-dimensional convex optimization problem.  We begin by expressing $\varphi^*$ as the solution to a high-dimensional convex optimization problem and then demonstrate how to reduce this problem to a more tractable two-dimensional one.

\Cref{eq:varphistar} defines $\varphi^*$ as the solution to a maximin problem that involves maximizing over distributions and minimizing over values of $t>0$.  The inner one-dimensional minimization over $t$ can be commuted with the outer infinite-dimensional maximization over $\boldsymbol{p}\in\mathcal{M}_{\boldsymbol{b},\mu}$.  Once the minimization and maximization operators have been commuted, the inner problem reduces to a finite-dimensional optimization problem.  This leads to a new expression for $\varphi^*$, as shown below.

\begin{theorem} \label{thm:tightchernoff} Fix $\boldsymbol{b}=(b_1,\ldots b_n)$, $\mu \in [0,\sum_{i=1}^n b_i]$, and $s \in [0,\sum_{i=1}^n b_i]$.  Let 
\begin{align*}
T=\left\{\tau\in \prod_{i=1}^m [0,b_i] : \sum_{i=1}^m \tau_i = \mu\right\}.
\end{align*}
Then, $\varphi^*_{\boldsymbol{b},\mu}(s)$ from \Cref{eq:varphistar} can be expressed as 
\begin{gather*}
\varphi^*_{\boldsymbol{b},\mu}(s)
= \min_{t\geq 0} \left(\max_{\tau \in T} \sum_{i=1}^m \log\left(1+\tau_i (\exp(b_i t)-1)/b_i \right) - ts\right).
\end{gather*}
\end{theorem}

The high-dimensional minimax problem from \Cref{thm:tightchernoff} can be reduced to a two-dimensional convex optimization problem.  

\begin{theorem} \label{thm:lowd} Fix $\boldsymbol{b}=(b_1,\ldots b_n)$, $\mu \in [0,\sum_{i=1}^n b_i]$, and $s \in [0,\sum_{i=1}^n b_i]$.  Let $\xi(b,t) = (\exp(bt)-1)/a$, 
\[
\tau^*_i(t,\lambda) = \min\left(\max\left(\frac{\xi(b_i,t)-\lambda}{\xi(b_i,t)\lambda},0\right),b_i\right),
\]
and
\[
g(t,\lambda;s) = \sum_i \log\left(1+\xi(b_i,t)\tau_i^*(t,\lambda)\right) + \lambda\left(\mu - \sum_i \tau_i^*(t,\lambda)\right)- ts.
\]
Then, $g$ is convex and $\varphi^*_{\boldsymbol{b},\mu}(s)$ from \Cref{eq:varphistar} can be expressed as 
\[
\varphi^*_{\boldsymbol{b},\mu}(s) = \min_{t,\lambda \geq 0} g(t,\lambda;s).
\]
Moreover, the mapping $t \mapsto \min_\lambda g(t,\lambda;s)$ is convex and the mapping $\lambda \mapsto g(t,\lambda)$ is convex for each $t \geq 0$.
\end{theorem}

This theorem suggests two methods for computing $\varphi^*$.  One method involves directly minimizing the two-dimensional convex function $t,\lambda \mapsto g(t,\lambda;s)$ and the other involves minimizing the one-dimensional convex function $t \mapsto \min_\lambda g(t,\lambda;s)$. Because each evaluation of the objective function of the latter requires solving a different one-dimensional convex optimization problem, in our simulations, we use the first method.

\section{Simulations}

\label{sec:simulations}

\begin{figure}
\begin{center}
\includegraphics[width=5in]{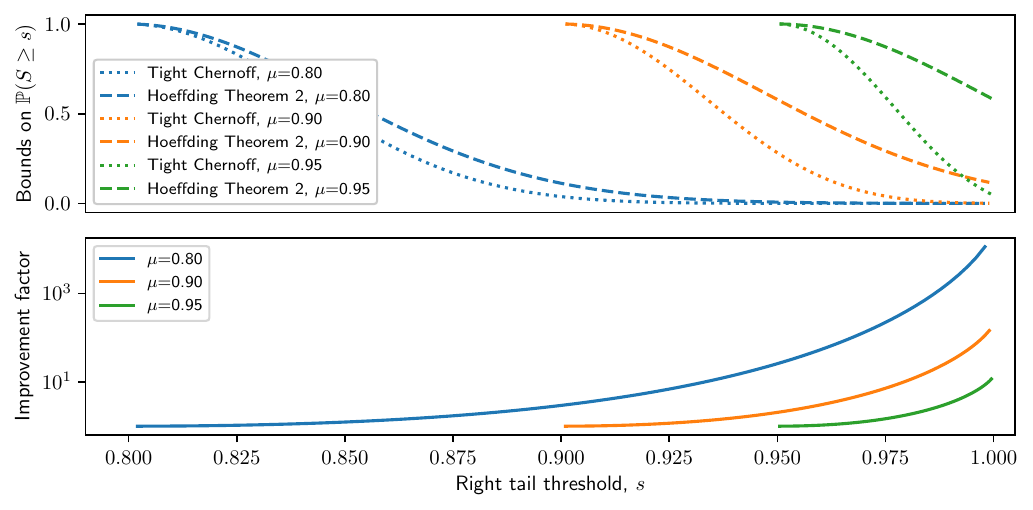}
  \caption{
  The new bounds are tighter than the general Hoeffding inequality.  We consider three choices of fixed expected value, $\mu \in \{0.8, 0.9, 0.95\}$.  We compare tail bounds based on the general Hoeffding inequality \citep[Theorem 2]{hoeffding_probability_1963} with the new bounds; the latter are calculated using \Cref{thm:lowd}. The top plot shows the bounds.  The bottom plot demonstrates the factor of improvement, showing ratios of the two bounds on a log scale.
  }\label{fig:tailboundcomparisons}
\end{center}
\end{figure}

We use simulations to assess the degree of improvement over the general Hoeffding inequality.  First, we sample a vector interval lengths, $\boldsymbol{b} \in \mathbb{R}^{200}$, uniformly at random from the set of lengths such that $B=\sum_{i=1}^{100} b_i=1$.  We define $S$ to be the sum of 200 independent random variables, $(X_1,\ldots,X_{200})$, where each $X_i$ is constrained to the interval $[0,b_i]$ for $i\in \{1,\ldots,200\}$.  

For any fixed value of $\mu$ and any value $s \geq \mu$, our new results provide an upper bound on $\mathbb{P}(S>s)$, assuming $\mathbb{E}[S]=\mu$.  The general Hoeffding inequality can also be used to produce a bound on $\mathbb{P}(S>s)$, offering a natural point of comparison.  These bounds on $\mathbb{P}(S>s)$ have practical implications, such as facilitating hypothesis testing for $\mathbb{E}[S]=\mu$ given an observed realization of $S$.

\Cref{fig:tailboundcomparisons} compares our new bounds with the general Hoeffding inequality for $\mu \in \{0.8, 0.9, 0.95\}$, showing how they vary as a function of $s$. In many cases, our bounds are orders of magnitude tighter than those derived from the general Hoeffding inequality, demonstrating their advantage in tail estimation.

\begin{appendix}
\section*{Proofs}\label{appn}

To prove \Cref{thm:tightchernoff}, it is convenient to view $\mathcal{M}_{\boldsymbol{b},\mu}$ from \Cref{eq:mdef} as a topological vector space.  For any scalar $b\geq 0$, let $\mathcal{M}_b$ denote the space of probability measures with support on $[0,b]$.  Endowing $\mathcal{M}_b$ with the weak-$*$ topology, we view $\mathcal{M}_b$ as a topological vector space.  Note that $\mathcal{M}_b$ is compact and convex.  For any values $\boldsymbol{b}\in ({\mathbb{R}^+})^m$ and fixed expected value $\mu$, let
\[
\mathcal{M}_{\boldsymbol{b},\mu} = \left\{\boldsymbol{p} \in \prod_i \mathcal{M}_{b_i}:\ \sum_i\mathbb{E}_{X_i \sim p_i}[X_i] \leq \mu\right\}.
\]
Here
$\prod_i \mathcal{M}_{b_i}$
denotes a direct sum of a topological vector spaces.  For example, if $\boldsymbol{p} \in \mathcal{M}_{\boldsymbol{b},\mu}$ then $\boldsymbol{p}=(p_1,\ldots, p_m)$ where each $p_i \in \mathcal{M}_{b_i}$.  Recall that addition is computed in the following manner for such direct sum spaces: for any $\boldsymbol{p},\boldsymbol{q} \in \mathcal{M}_{\boldsymbol{b},\mu}$, $(\boldsymbol{p}+\boldsymbol{q})=((p_1+q_1),\ldots, (p_m+q_m))$.  This seemingly excessive formality is necessary to clarify that $\mathcal{M}_{\boldsymbol{b},\mu}$ is, in fact, convex.

\begin{theorem*}[Restatement of Theorem 3.1] Fix $\boldsymbol{b}=(b_1,\ldots b_n)$, $\mu \in [0,\sum_{i=1}^n b_i]$, and $s \in [0,\sum_{i=1}^n b_i]$.  Let $T=\{\tau\in \prod_{i=1}^m [0,b_i]: \sum_{i=1}^m \tau_i = \mu\}$.   Then $\varphi^*_{\boldsymbol{b},\mu}(s)$ from \Cref{eq:varphistar} can be expressed as 
\begin{gather*}
\varphi^*_{\boldsymbol{b},\mu}(s)
= \min_{t\geq 0} \left(\max_{\tau \in T} \sum_{i=1}^m \log\left(1+\tau_i (\exp(b_i t)-1)/b_i \right) - ts\right).
\end{gather*}
\end{theorem*}

\begin{proof}

Fix $s$. Let $\xi(b,t)=(\exp(bt)-1)/a$.  Let 
\[
f(t,\boldsymbol{p})=\sum_{i}\log \mathbb{E}_{B_{i}\sim p_{i}}[\exp(tB_{i})]-ts.
\]
Our object of interest is may be given as $\min_{t}\max_{\boldsymbol{p}}f(t,\boldsymbol{p})$.
Our first step is to reverse the order of the minimization and maximization
by applying Sion's minimax theorem \citep{komiya1988elementary}. We observe the following.
\begin{itemize}
\item $f$ is convex with respect to $t$, as it is the sum of cumulant
generating functions. 
\item $f$ is concave with respect to $\boldsymbol{p}$, due to Jensen's
inequality.
\item $f$ is continuous with respect both $\boldsymbol{p}$ and $t$.
\item $\mathcal{M}_{\boldsymbol{b},\mu}$ is both convex and compact.
\item $[0,\infty)$ is convex (though not compact).
\end{itemize}
Sion's minimax theorem thus shows that
\[
\min_{t\geq0}\max_{\boldsymbol{p}\in\mathcal{M}_{\boldsymbol{b},\mu}}f(t,\boldsymbol{p})=\inf_{t\geq0}\max_{\boldsymbol{p}\in\mathcal{M}_{\boldsymbol{b},\mu}}f(t,\boldsymbol{p}).
\]
For any fixed $t$, the maximization of $f$ may be reduced to a finite-dimensional
problem. First, we rewrite the maximization problem by introducing
auxiliary variables, as follows:
\begin{align*}
\max_{\boldsymbol{p},\tau}\qquad & f(t,\boldsymbol{p})\\
\mathrm{s.t.}\qquad & \tau_{i}\in[0,b_i]\\
 & p_{i}\in\mathcal{M}_{b_i}\\
 & \tau_{i}=\sum_{i}\mathbb{E}_{B_{i}\sim p_{i}}\left[B_{i}\right]\\
 & \sum_{i}\tau_{i} \leq \mu
\end{align*}
For any fixed $\tau$, Hoeffding \citep[proof of Theorem 1]{hoeffding_probability_1963} showed that the cumulant generating
functions can be maximized by setting each $p_{i}$ to be a scaled
Bernoulli random variable, namely
\[
p_{i}=\delta_{0}\left(1-\frac{\tau_{i}}{b_i}\right)+\delta_{b_i}\frac{\tau_{i}}{b_i}
\]
where $\delta_{x}$ represents the point mass at $x$. Under this
distribution, 
\[
\log \mathbb{E}_{p_{i}}\left[\exp\left(tB_{i}\right)\right]=\log\left(1+\xi(b_i,t)\tau_{i}\right).
\]
Note that this is monotone increasing in $\tau$, so the constraint $\sum_i \tau_i \leq \mu$ may be replaced with the constraint $\sum_i \tau_i=\mu$ without changing the result.  Our problem thus reduces to the desired form.

\end{proof}

\begin{theorem*}[Restatement of Theorem 3.2] Fix $\boldsymbol{b}=(b_1,\ldots b_n)$, $\mu \in [0,\sum_{i=1}^n b_i]$, and $s \in [0,\sum_{i=1}^n b_i]$.  Let $\xi(b,t) = (\exp(bt)-1)/a$, 
\[
\tau^*_i(t,\lambda) = \min\left(\max\left(\frac{\xi(b_i,t)-\lambda}{\xi(b_i,t)\lambda},0\right),b_i\right),
\]
and
\[
g(t,\lambda;s) = \sum_i \log\left(1+\xi(b_i,t)\tau_i^*(t,\lambda)\right) + \lambda\left(\mu - \sum_i \tau_i^*(t,\lambda)\right)- ts.
\]
Then, $g$ is convex and $\varphi^*_{\boldsymbol{b},\mu}(s)$ from \Cref{eq:varphistar} can be expressed as 
\[
\varphi^*_{\boldsymbol{b},\mu}(s) = \min_{t,\lambda \geq 0} g(t,\lambda;s).
\]
Moreover, the mapping $t \mapsto \min_\lambda g(t,\lambda;s)$ is convex and the mapping $\lambda \mapsto g(t,\lambda)$ is convex for each $t \geq 0$.
\end{theorem*}

\begin{proof}
Let $T=\{\tau\in \prod_{i=1}^m [0,b_i]: \sum_{i=1}^m \tau_i = \mu\}$.  Recall from \Cref{thm:tightchernoff} that 
\begin{gather*}
\varphi^*_{\boldsymbol{b},\mu}(s)
= \min_{t\geq 0} \left(\max_{\tau \in T} \sum_{i=1}^m \log\left(1+\tau_i (\exp(b_i t)-1)/b_i \right) - ts\right).
\end{gather*}
For any fixed $t\geq0$, we begin by rewriting the inner maximization problem
in terms of 
\[
f(\tau;t)=\begin{cases}
\sum_{i}\log\left(1+\xi(b_i,t)\tau_{i}\right)-ts & \mathrm{if}\ \tau\in\prod_{i}[0,b_i]\\
-\infty & \mathrm{else}.
\end{cases}
\]
The inner maximization problem is thus equivalent to
\[
\max_{\tau:\ \sum_{i}\tau_{i}=\mu}f(\tau;t).
\]
The associated Lagrangian function is given by 
\[
\mathcal{L}(\tau,\lambda;t)=\sum_{i}\log\left(1+\xi(b_i,t)\tau_{i}\right)+\lambda\left(\mu-\sum_{i}\tau_{i}\right)-ts.
\]
For any fixed $t\geq0$, observe that $\xi(b_i,t)$ is positive
and thus $\tau\mapsto\mathcal{L}\left(\tau,\lambda;t\right)$ is concave.
The argument maximizing $\tau\mapsto\mathcal{L}(\tau,\lambda;t)$
is given by $\tau_{i}^{*}(t,\lambda)$. Thus $g(t,\lambda)=\max_{\tau}\mathcal{L}(\tau,\lambda;t)$,
the Lagrangian dual for a convex optimization problem with convex
constraints. There is at least one feasible point, namely $\tau_{i}=b_i\mu/\sum_{j}a_{j}$,
so strong duality holds and 
\[
\max_{\substack{\tau\in\prod_{i}[0,b_i]\\
\sum_{i}\tau_{i}=\mu
}
}\sum_{i}\log\left(1+\xi(b_i,t)\tau_{i}\right)-ts=\min_{\lambda}g(t,\lambda)
\]
as desired.

We now demonstrate that $g$ is convex. First, observe that $\mathcal{L}$
is affine in $\lambda$ and convex in $t$. Thus $g$ is a pointwise
maximum of a family of convex functions: it is convex. 

Finally, we demonstrate that $t\mapsto\min_{\lambda}g(t,\lambda)$
is convex. Observe that, for any fixed $\tau_{i}\geq0$, the mapping
\[
t\mapsto\sum_{i}\log\left(1+\xi(b_i,t)\tau_{i}\right)-ts
\]
is convex in $t$. Thus $\min_{\lambda}g(t,\lambda)$ is also the
pointwise maximum of a family of convex functions: it is also convex.
\end{proof}  

\end{appendix}

\bigskip
\begin{center}
{\large\bf SUPPLEMENTARY MATERIAL}
\end{center}

\begin{description}

\item[Demonstrating the new bound:] Code for solving the two-dimensional convex optimization problem involved in the new bound, together with demonstrations of its use.  (Jupyter notebook)

\end{description}

\bibliographystyle{agsm}
\bibliography{bibliography}       

\end{document}